\documentclass[leqno]{article}
\usepackage{epsfig}
\input epsf
\usepackage{latexsym}
\usepackage{amsmath}

\begin{document}

\title{Helicoidal surfaces with constant anisotropic mean curvature}
\author{By C{\footnotesize HAD} K{\footnotesize UHNS}
and B{\footnotesize ENNETT} P{\footnotesize ALMER}}
\date{}
\maketitle
\newtheorem{theorem}{Theorem}[section]
\newtheorem{cor}{Corollary}[section]
\newtheorem{prop}{Proposition}[section]
\newtheorem{lemma}{Lemma}[section]
\newtheorem{condition}{Condition}[section]
\newtheorem{example}{Example}[section]

\newtheorem{definition}{Definition}[section]
\newtheorem{remark}{Remark}[section]
\newtheorem{conjecture}{Conjecture}[section]
\newtheorem{claim}{Claim}[section]
\newtheorem{question}{Question}[section]
\newcommand{\rf}[1]{\mbox{(\ref{#1})}}

\renewcommand{\thefootnote}{\fnsymbol{footnote}}

\abstract{We study surfaces with constant anisotropic mean curvature which are invariant under a helicoidal motion. For functionals with axially symmetric Wulff shapes, we generalize the recently developed twizzler representation (\cite{P}) to the anisotropic case and show how all helicoidal constant anisotropic mean curvature surfaces can be obtained by quadratures.}\\[5mm]

Helicoidal symmetry occurs in a wide range of physical phenomena. It occurs frequently in biological systems \cite{BF} due to the fact  that it  arises
from a fundamental self-organizational principle that a regular assembly of identical objects has helical symmetry.  On the microscopic scale, instances of helicoidal symmetry include the orientation of molecules of cholesteric liquid crystal
and certain twist grain boundaries of diblock polymers. 

Anisotropic surface energies occur at interfaces between immiscible materials when at least one of them is in an ordered phase. The
simplest example is the  free energy
\begin{equation}
\label{F} {\cal F}=\int_\Sigma \gamma (\nu)\:d\Sigma\:,\end{equation}
where $\gamma$ is the anisotropic energy density which gives the unit energy per unit area of a surface element having unit normal vector $\nu$.
Such energies were first applied to study the free surface energy of crystals. Wulff stated that the equilibrium shape of a crystal could be obtained by minimizing a specific anisotropic surface energy subject to a volume constraint.

We consider a surface  given as a smooth, oriented immersion
$X:\Sigma \rightarrow {\bf R}^3$ with Gauss map
$\nu:\Sigma \rightarrow S^2$. For a smooth variation of $X$,
$X_\epsilon =X+\epsilon {\dot X}+...$, we have the first variation
formula,
$$\delta {\cal F}=:- \int_\Sigma \Lambda {\dot X}\cdot \nu\:d\Sigma +{\rm boundary\:\:terms}\:.$$
This formula defines the {\it anisotropic mean curvature}
 $\Lambda$. If $V$ denotes the algebraic volume enclosed by the
surface:
$$V=\frac{1}{3}\int_\Sigma X\cdot \nu\:d\Sigma\:,$$
then a well known formula for the variation of $V$ is
$$\delta V=\int_\Sigma {\dot X}\cdot \nu\:d\Sigma\:.$$
Therefore $\Lambda\equiv $ constant characterizes the volume constrained equilibria of ${\cal F}$.

Following \cite{G}, we now give a way to locally calculate of the anisotropic mean curvature.  If the function $\gamma$ is sufficiently smooth,
let ${\tilde \gamma}$ denote the positive homogeneous degree one extension of $\gamma$, i.e. ${\tilde \gamma }(Y):=|Y|\gamma(Y/|Y|)$ for $Y\in {\bf R}^3\setminus \{0\}$.
The {\it Cahn-Hoffman} field is defined by
$$\xi_p:=(\nabla {\tilde \gamma})_{\nu(p)}\:.$$
This field, which is always transversal to the surface, can be thought of as an anisotropic normal
field to to the surface.
The anisotropic mean curvature is then given by
\begin{equation}
\label{lambda}
\Lambda=-{\rm Div}\: \xi\:.\end{equation}
Here the divergence can be computed on the surface or in three dimensional space if $\xi$ is first smoothly extended to a field near the surface.

There is a canonical surface associated with the anisotropic energy density $\gamma$ called the {\it Wulff shape}, which can be defined by
\begin{equation}
\label{W} W=\partial  \bigcap_{n\in S^2} \{ Y\cdot n\le \gamma(n)\}\:.\end{equation}
(Some authors define the Wulff shape as the intersection itself.)
As the boundary of an intersection of half-spaces, $W$ is convex and it will be assumed in this paper that $W$ is smooth and has uniformly positive curvature $K_W>0$. With this assumption,
the equation \rf{lambda}, with the right hand side being any function of the space variables, is elliptic.  A fundamental result, known as Wulff's Theorem, roughly  states that $W$ is the absolute minimizer of the free energy ${\cal F}$ among all closed surfaces enclosing the same volume, thus $W$ solves the isoperimetric
problem for the anisotropic energy functional ${\cal F}$. From this it follows that $W$ has constant anisotropic mean curvature.

 It should be noted that any closed convex surface can be realized as the Wulff shape for some functional. If $W$ is any convex surface, then its Gauss map is a diffeomorphism and its Gauss map $N:W\rightarrow S^2$ is a diffeomorphism. The functional with anisotropic density function
given by $\gamma:=(r \cdot N)\circ N^{-1}$, where $r$ is the position vector on $W$, then has Wulff shape $W$. This is a useful construction since it is sometimes more convenient to specify the Wulff shape instead of producing a formula for the density function. 

 The purpose of this paper is to discuss equilibrium surfaces for a volume constrained anisotropic energy of the form \rf{F} which are invariant under a helicoidal motion. The isotropic case of constant mean curvature helicoidal surfaces ($\gamma \equiv 1$), has been discussed in \cite{dCD}. \cite{R}, \cite{S}, \cite{P}. The helicoidal surfaces with constant mean curvature arise as isometric deformations of Delaunay surfaces. In this deformation the principal curvatures are preserved. The best known example is the isometric deformation between the catenoid and the helicoid.

In this paper  show a certain universal property of the classical helicoid. It has zero anisotropic mean curvature for every rotationally symmetric anisotropic energy of the type discussed above. We then derive the equilibrium equations for the constant anisotropic mean curvature surfaces which are invariant under a helicoidal motion; the so called ``twizzlers''. We do this first for surfaces in non parametric form. Following this, we generalize a representation formula developed by Perdomo \cite{P} for the isotropic (CMC) case. Finally we generalize the conjugacy relation between the catenoid and the helicoid to a wide range of anisotropic functionals.

\section{Euler-Lagrange equation in non parametric form}
Let $W\subset {\bf R}^3$ be a smooth convex surface. There exists an embedding $\chi$ of $S^2$ into ${\bf R}^3$ such that $\chi(S^2)=W$ and $(\chi|_W)^{-1}$ is the
Gauss map of $W$.

Suppose $\Sigma$ is given as the graph $z=z(x,y)$ over a domain $D$ in the plane. The normal map $\nu$  will take
values in a hemisphere and we can choose the orientation so that it is the upper hemisphere. 
If we consider the composition
$$\Sigma \xrightarrow{\nu} S^2 \xrightarrow{\chi}{\bf R}^3\:,$$
then $\xi (\Sigma)$ will lie in a part of $W$ for which the normal map $\chi^{-1}$ lies
in the upper hemisphere. Then this part of $W$ can be represented as a graph
$$\xi_1, \xi_2 \mapsto (\xi_1, \xi_2,v(\xi_1, \xi_2)\:.$$
If the curvature of $W$ is strictly positive, then we have
$$v_{\xi_1 \xi_1}v_{\xi_2 \xi_2}-v_{\xi_1 \xi_2}^2>0
$$
and it is possible to write
\begin{equation}
\label{inv}
\xi_1=\xi_1(v_{\xi_1}, v_{\xi_2}),\qquad \xi_2=\xi_2(v_{\xi_1}, v_{\xi_2})\:.
\end{equation}
By composing $v(\xi_1, \xi_2)$ with the transformation (\ref{inv}) we obtain a function
$$V(v_{\xi_1}, v_{\xi_2}):=v({\xi_1}(v_{\xi_1}, v_{\xi_2}),{\xi_2}(v_{\xi_1}, v_{\xi_2}))\:.$$

At points where the normals on $\Sigma$ and $W$ agree, we have
$$\frac{(-z_x, -z_y, 1)}{\sqrt{1+z_x^2+z_y^2}}=\frac{(-v_{\xi_1}, -v_{\xi_2}, 1)}{\sqrt{1+v_{\xi_1}^2+v_{\xi_2}^2}}\:$$
from which there follows
\begin{equation}
\label{par}
z_x=v_{\xi_1}\:,\qquad z_y=v_{\xi_2}\:.
\end{equation}
Let $\gamma$ denote the support function of $W$. Then
$$\gamma=\langle ({\xi_1}, {\xi_2}, v)  ,\frac{(-v_{\xi_1}, -v_{\xi_2}, 1)}{\sqrt{1+v_{\xi_1}^2+v_{\xi_2}^2}}\rangle=
\frac{v-{\xi_1} v_{\xi_1}-{\xi_2} v_{\xi_2}}{\sqrt{1+v_{\xi_1}^2+v_{\xi_2}^2}}\:.$$
It then follows from (\ref{par}), that the energy is given by
$${\cal F}=\int_D (V-\xi_1z_x-\xi_2 z_y)\:dxdy\:,$$
where $\xi_1=\xi_1(z_x, z_y)$, $\xi_2=\xi_2(z_x, z_y)$ and $V=V(z_x, z_y)$ are defined above.

We want to compute the first variation.  Replacing $z$ by $z+\epsilon {\dot z}$ and taking the
derivative of the integrand in ${\cal F}$ with respect to $\epsilon$, gives
\begin{eqnarray*}
\partial_\epsilon (V-{\xi_1} z_x-{\xi_2} z_y)_{\epsilon=0}&=& V_{v_{\xi_1}}{\dot z_x}+V_{v_{\xi_2}}{\dot z_y}
-{\xi_1}_{v_{\xi_1}}z_x{\dot z_x}-{\xi_1}_{v_{\xi_2}}z_x{\dot z_y}-{\xi_1} {\dot z_x}-{\xi_2} {\dot z_y}\\
&=& {\dot z_x}[V_{v_{\xi_1}}-z_x {\xi_1}_{v_{\xi_1}}-z_y{\xi_1}_{v_{\xi_1}}-{\xi_1}]+{\dot z_y}[V_{v_{\xi_2}}-z_x {\xi_1}_{v_{\xi_2}}-z_y{\xi_2}_{v_{\xi_2}}-{\xi_2}]
\end{eqnarray*}
However, the chain rule , \rf{par} and the definition of $V$, gives
$$V_{v_{\xi_1}}=v_{\xi_1} {\xi_1}_{v_{\xi_1}}+v_{\xi_2} {\xi_2} _{v_{\xi_1}}=z_x  {\xi_1}_{v_{\xi_1}}+z_y{\xi_2} _{v_{\xi_1}}\:$$
and
$$V_{v_{\xi_2}}=v_{\xi_1} {\xi_1}_{v_{\xi_2}}+v_{\xi_2} {\xi_2} _{v_{\xi_2}}=z_x  {\xi_1}_{v_{\xi_2}}+z_y{\xi_2} _{v_{\xi_2}}\:$$

Using this above, gives
$$\partial_\epsilon (V-{\xi_1} z_x-{\xi_2} z_y)_{\epsilon=0}= {\dot z_x}[-{\xi_1}]+{\dot z_y}[-{\xi_2}]\:,$$
so that, assuming ${\dot z}$ has compact support in $D$,
\begin{equation}
\label{var}
\delta{\cal F}=\int_D {\dot  z}({\xi_1}_x+{\xi_2}_y)\:dxdy\:.
\end{equation}
i.e. $\Lambda \equiv 0$ is equivalent to the equation
\begin{equation}
\label{d2}
{\rm Div}_0 ({\xi_1}(z_x,z_y), {\xi_2}(z_x, z_y))=0\:,
\end{equation}
where ${\rm Div}_0$ denotes the two dimensional divergence.

It is clear from \rf{W} that $\gamma$ is the support function of the Wulff shape $W$. The classical representation of a convex surface
by its support function, known as the {\it tangential representation} \cite{E}, then gives:
$$\xi=D\gamma+\gamma \nu\:.$$
In the case where $\gamma=\gamma(\nu_3)$, i.e. when $W$ is axially symmetric, we obtain $\xi=\gamma'(\nu_3)E_3^T+\gamma \nu$. Letting $1/\mu_2:=\gamma(\nu_3)-\nu_3\gamma'(\nu_3)$, we can have
\begin{eqnarray*}\xi&=&D\gamma +\gamma \nu_3\\
&=& \gamma'(\nu_3)E_3^T+\gamma(\nu_3)\nu\\
&=&\gamma'(\nu_3)(E_3-\nu_3 \nu)+\gamma(\nu_3)\nu\\
&=&\frac{1}{\mu_2} \nu+\gamma'(\nu_3)E_3\:.\end{eqnarray*}
In particular $(\xi_1, \xi_2)=(1/\mu_2)(\nu_1,\nu_2)$.  Since the normal $\nu$ is given
by $\nu=\nu_3(-z_x, -z_y, 1)$, we arrive at the Euler-Lagrange equation
\begin{equation}
\label{ELgraph}
{\rm Div}_0\biggl[ \frac{\nu_3}{\mu_2} (z_x, z_y)\biggr]=\Lambda \:(\equiv {\rm constant})\:\end{equation}
for the volume constrained variational problem. Here $\nu_3=(1+|\nabla z|^2)^{-1/2}$.
\section{Helicoidal solutions}
We will seek a solution of \rf{ELgraph} of the form $z=g(r)+\Lambda \theta$, where $\lambda$ is a real constant. In this case, the surface is given
by $X=(re^{i\theta}, g(r)+\lambda \theta)$, where we have replaced the first two coordinates in ${\bf R}^3$ with a complex coordinate. We then
obtain $|\nabla z|^2=g_r^2+\lambda^2/r^2$, so that $\nu_3$ and $\mu_2(\nu_3)$ are independent of $\theta$. Equation \rf{ELgraph} can be expressed
$$\Lambda r\:dr\:d\theta =d* \frac{\nu_3}{\mu_2}(rg_r dr+\lambda d\theta)=d\:\frac{\nu_3}{\mu_2}(rg_r \:d\theta -\frac{\lambda}{r}\:dr)\:,$$
from which it follows easily that $(\nu_3rg_r)/\mu_2)_r =\Lambda r$ holds. Integrating, we obtain the first integral
\begin{equation}
\label{int} \frac{\nu_3rg_r}{\mu_2(\nu_3)}-\frac{\Lambda r^2}{2}=C\:,\end{equation}
where $C$ is a constant.  Note that in \rf{int},  the helicity of the surface  is built into the dependence of ${\nu_3}$ in $\lambda$. When $\lambda=0$, we can recover the representation
of anisotropic Delaunay surfaces which was found in \cite{KP2005}.

We obtain from \rf{int}:
\begin{prop} For any axially symmetric anisotropic energy density $\gamma=\gamma(\nu_3)$, the usual helicoids given by $z=\lambda \theta +C$, $\lambda, C\in {\bf R}$ have zero anisotropic mean curvature.
\end{prop}

The Cahn-Hoffman field defines a map $\xi:\Sigma \rightarrow W$ which can be considered as a type of anisotropic Gauss map. In the special case
that both $W$ and $\Sigma$ are axially symmetric and $\Sigma$ has constant anisotropic mean curvature, this map is harmonic \cite{KP2008}. We will use the example of the helicoid to show that, in general, constancy of $\Lambda$ is not enough to insure harmonicity of $\xi$.

We assume that $W$ is any axially symmetric Wulff shape and we let $\Sigma$ be a helicoid given as $X=(re^{i\theta}, C+\lambda\theta)$. The harmonicity of $\xi$ is equivalent to
 $(\Delta \xi)^T=0$, where $\Delta$ is the Laplacian and the superscript $T$ denotes the tangential component. Since $W$ is axially symmetric, it follows from the remarks above that
 $$\xi =\frac{1}{\mu_2} \nu +\gamma'(\nu_3)E_3\:,$$
 So
 $$\Delta \xi =(\Delta \frac{1}{\mu_2})\nu +2d\nu(\nabla \frac{1}{\mu_2})+\frac{1}{\mu_2}\Delta \nu+\Delta \gamma'(\nu_3)E_3\:.$$
 By a well known formula  $(\Delta \nu)^T=-2\nabla H =0$ since the helicoid is a minimal surface. We get,
\begin{eqnarray*}(\Delta \xi)^T&=&2d\nu(\nabla \frac{1}{\mu_2})+(\Delta \gamma'(\nu_3) E^T_3\\
&=& ( \frac{1}{\mu_2})'(\nu_2)2d\nu(\nabla \nu_3)+\gamma''(\nu_3)\Delta \nu_3 +\gamma '''(\nu_3)|\nabla \nu_3|^2\:E_3^T\\
&=&\biggl( (-\nu_3\gamma ''(\nu_3))(-K)+\gamma'''(\nu_3)|\nu_3|^2\biggr)\:E_3^T,
\end{eqnarray*}
where we have used that $\nabla \nu_3=-K E_3^T$ on a minimal surface. A calculation shows that $K=-\lambda^2(r^2+\lambda^2)^{-2}$ on $\Sigma$. Using that the the normal to $X$ is
given by $\nu=(r^2+\lambda^2)^{-1/2}(-i\lambda e^{i\theta}, r)$, we arrive at the following.
\begin{prop}
The Cahn-Hoffman field of  a helicoid defines a harmonic map into the Wulff shape if and only if
the anisotropic energy density $\gamma$ satisfies the differential equation
$$\gamma'''=\frac{4}{\lambda^2}\nu_3(1-\nu_3^2)\gamma''\:.$$
\end{prop}

Although the Cahn-Hoffman map of the helicoid is not in general harmonic, it is a critical point of the energy if one allows the metric tensor to vary as the immersion varies, i. e. it is a critical point of the action given locally by
$$X\mapsto \int g^{ij}\xi_{,i} \xi_{,j} \sqrt{g} \:d^2x\:,$$
where $(g_{ij}$ are the components of the induced metric from $X$. We refer the reader to \cite{pa}
for details.

We will now consider a particular functional for which the integration of the equation \rf{ELgraph} is particularly easy.
Using $\nu_3=(1+|\nabla z|^2)^{-1/2}$ and $d\Sigma =\nu_3^{-1}\:d^2x$,  the Dirichlet integral can be expressed
$$D[z]:=\frac{1}{2}\int_\Omega |\nabla z|^2\:d^2x=
\int_\Sigma \frac{1}{\nu_2}-\nu_3\:d\Sigma\:.$$
This means that the functional with density $\gamma:=\nu_3^{-1}-\nu_3$
possesses critical points which are graphs of harmonic functions.
In this case, $1/\mu_2=1/\nu_3$ and so \rf{int} reduces to
$rg_r-\Lambda r^2/2=C$. Integration yields,
$$g=\frac{\Lambda r^2}{4}+C\log r+C_1\:.$$

The Wulff shape corresponding to the density $\gamma$ is the elliptic parabloid $\xi_1^2+\xi_2^2=-2\xi_3$.
Some examples are shown below.
\begin{figure}[h]
\begin{minipage}[c]{.45\textwidth}
    \begin{center}

       \includegraphics[width=55mm,height=55mm,angle=0]{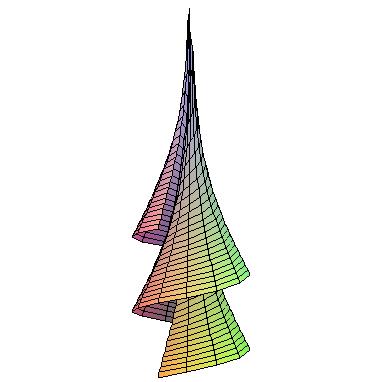}

     % \label{c1}
    \end{center}
\end{minipage}
  %\hfill
  \begin{minipage}[c]{.45\textwidth}
    \begin{center}

      \includegraphics[width=55mm,height=55mm,angle=0]{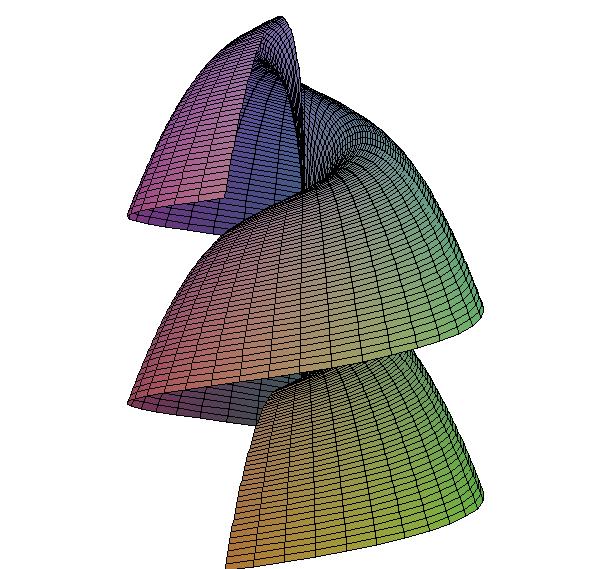}

     % \label{c2}
    \end{center}
   \end{minipage}
 % \caption{Wulff shape (left) and anisotropic catenoid (right)}
 % \label{dels1}
  %\end{figure}

  %\hfill
   \begin{minipage}[c]{.45\textwidth}
  \begin{center}
        \includegraphics[width=55mm,height=55mm,angle=0]{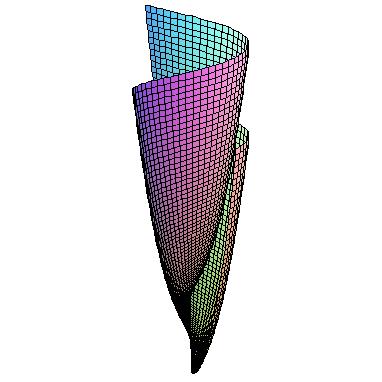}

     % \label{c2}
    \end{center}
   \end{minipage}
    \begin{minipage}[c]{.45\textwidth}
  \begin{center}
        \includegraphics[width=55mm,height=55mm,angle=0]{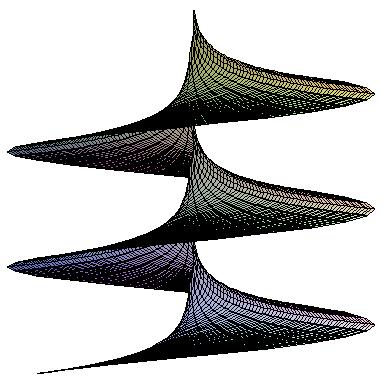}

     % \label{c2}
    \end{center}
   \end{minipage}
   \end{figure}
\section{Twizzler representation}
In this section we develop a representation for helicoidal surfaces with constant anisotropic mean curvature. Our treatment is based on \cite{P} in which the author derives a representation formula in the isotropic (i.e. constant mean curvature) case.

We will write a helicoidal surface in the form
\begin{equation}
\label{tw} X=([x(s)+iy(s)]e^{-i\omega \vartheta}, \vartheta +C)\:.\end{equation}
Here $\omega$ and $C$ are real constants and $(x(s), y(s))$ is a plane curve parameterized by
arc length. We will refer to this curve as the {\it generating curve} of $X$ since the surface is the orbit of this curve under the helicoidal action. Above, we have identified the $x_1x_2$ plane with the complex plane ${\bf C}$. The surface given by \rf{tw} is clearly invariant under the ${\bf R}$ action:
$$(\:t\:,\:(x+iy,z)\:)\mapsto (e^{i\omega t}(x+iy), z-t)\:.$$

For the given curve $(x(s), y(s))$, we define
\begin{equation}
\label{etas}\eta_1(s):=\frac{-1}{2}\partial_s [(x(s))^2+(y(s))^2]\:,\qquad \eta_2(s):=xy_s-yx_s\:.\end{equation}
Perdomo, \cite{P}, refers to the curve $(\eta_1(s), \eta_2(s))$ as the {\it treadmill sled} of the curve $(x(s), y(s))$. It is, acording to \cite{P}, the trace of the origin when the curve $(x(s), y(s))$
rolls without slipping on a treadmill located at the origin aligned along the x-axis.
\begin{theorem}
\label{T1}

The surface $X$ given by \rf{tw} has constant anisotropic mean curvature $\Lambda$ if and
only if the following equation holds
\begin{equation}
\label{H} \Lambda (\eta_1^2+\eta^2_2) +\frac{2\eta_2}{\mu_2(\frac{\omega\eta_1}{\sqrt{1+\omega^2\eta_1^2}})\: \sqrt{1+\omega^2\eta_1^2}}+ A\equiv 0\:,
\end{equation}
where $A$ is a real constant.
\end{theorem}

\begin{prop} The generating curve $(x,y)$ can be recovered from the curve in \rf{H} by the formula
\begin{equation}
\label{tc}
x+iy=-(\eta_1+i\eta_2)\exp\biggl( -i\int\frac{d\eta_2}{\eta_1}\biggr)\:.\end{equation}
\end{prop}
{\it Proof.}    Recall that $s$ is the arc length parameter of the curve $(x,y)$. If we write $x_s+iy_s=:\exp (i\phi(s))$,
then
\begin{equation}
d\phi =-\kappa ds\:,\end{equation}
holds where $\kappa$ denotes the curvature of the curve $(x,y)$. On the other hand, it is easy to check that
\begin{equation}
\frac{d\eta_2 }{\eta_1}=\kappa ds=-d\phi \:,\end{equation}
holds. Combining, we get
$$\exp\biggl( -i\int\frac{d\eta_2}{\eta_1}\biggr)=e^{i\phi}=x_s+iy_s \:.$$
Finally, a simple calculation using the definitions of $\eta_1$ and $\eta_2$ gives
\begin{equation}
\label{xiy}
x+iy= -(\eta_1+i\eta_2)(x_s+iy_s)\end{equation}
 which gives \rf{tc}. {\bf q.e.d}\\[4mm]
{\bf Remark} Theorem  \rf{T1} allows the construction of all helicoidal surfaces with $\Lambda=$ constant as follows:
\begin{itemize}
\item Regard \rf{H} is a quadratic in $\eta_2$, the equation can be solved for $\eta_2=\eta_2(\eta_1, \omega, \Lambda, A)$.
\item Generate the twizzler curve $x(s)+iy(s)$ using
\rf{tc}.
\item Recover the immersion $X$ using \rf{tw}. \\[2mm]
\end{itemize}
Some examples are shown in Figures 1 and 2 which are based on the Rapini-Papoular functional $\gamma=1+e \nu_3^2$, $e\in {\bf R}$.

\begin{lemma} Let $X:\Sigma \rightarrow {\bf R}^3$ be a helicoidal surface
given by \rf{tw} and assume that \rf{H} holds. If there exists an open set
$U\subset \Sigma$ such that $\nu_3\equiv 0$ on $U$, then all of
$X(\Sigma)$ is contained in a circular cylinder.
\end{lemma}
{\it Proof.} A calculation using \rf{tw} shows that 
\begin{equation}
\label{nu3}
\nu_3=\frac{\omega \eta_1}{\sqrt{1+\omega^2 \eta_1^2}}\:,\end{equation} so  $\nu_3\equiv 0$ on $U$ implies that $\eta_1\equiv 0$. It follows from \rf{etas} that locally $(x(s))^2+(y(s))^2$
is identically a non zero constant $=:R$. It then easily follows that $\eta_2=R$
and so an open set in the surface $U_1\subset U$ is contained in a vertical circular cylinder of radius $R$. 

The anisotropic mean curvature of the cylinder is given by $\Lambda =-1/(R\mu_2(0))$, (see \cite{KP2005}). We then obtain from \rf{tw}, 
$$-A=\Lambda \eta_2^2+\frac{2\eta_2}{\mu_3(0)}=
\frac{-R}{\mu_2(0)}+\frac{2R}{\mu_2(0)}=\frac{R}{\mu_2(0)}\:.$$
Therefore, on $U$ there holds
\begin{equation}
\label{LA}
 \Lambda A =\frac{1}{(\mu_2(0))^2}\:.\end{equation}
We will use this to show that the entire surface is contained in the cylinder of radius $R$.

Regarding \rf{H} as a quadratic equation for $\eta_2$, one sees that the discriminant is 
$$\frac{4}{(\mu_2(\nu_3))^2(1+\omega^2\eta_1^2)}-4\Lambda^2\eta_1^2-4\Lambda   A\ge 0$$
since $\eta_2$ is real. Using \rf{nu3} and \rf{LA}, one sees that the previous inequality is the same as
\begin{equation}
\label{A} \frac{1-\nu_3^2}{(\mu_2(\nu_3))^2}-\Lambda_1^2\eta_1^2 \ge \frac{1}{(\mu_2(0))^2}\:.\end{equation} 

Recall from \cite{KP2005} that the principal curvatures, $\mu_i$, $i=1,2$, of the Wulff shape $W$ with respect to the inward pointing normal are given by 
$$\frac{1}{\mu_2}=\gamma -\nu_3\gamma'(\nu_3)\:,\quad 
\frac{1}{\mu_1}=(1-\nu_3^2)\gamma''(\nu_3)+\frac{1}{\mu_2}\:.$$
By the convexity condition, $\mu_i>0$ holds on $W$.
Differentiation shows
\begin{eqnarray*}
\frac{d}{d \nu_3}\biggl(\frac{1-\nu_3^2}{(\mu_2(\nu_3))^2}\biggr)&=&\frac{d}{d \nu_3}\biggl((1-\nu_3^2)(\gamma -\nu_3\gamma')^2\biggr)\\
&=& -2\nu_3(\gamma -\nu_3\gamma')^2+(1-\nu_2)^22(\gamma -\nu_3\gamma')(-\nu_3\gamma'')\\
&=&-2\nu_3(\gamma -\nu_3\gamma')\bigl((1-\nu_3^2)\gamma''(\nu_3)+\gamma-\nu_3\gamma'\bigr)\\
&=&\frac{-2\nu_3}{\mu_1\mu_2}\:.\end{eqnarray*}
It follows that the derivative is negative for $\nu_3>0$ and positive for $\nu_3<0$, 
so $(1-\nu_3^2)/\mu_2(\nu_2)$ has a maximum at $\nu_3=0$. It then follows from 
\rf{A} that $\nu_3\equiv 0$ and $\eta_1\equiv 0$ holds on $\Sigma$. {\bf q.e.d}\\[2mm]
\begin{lemma}  Let $X:\Sigma \rightarrow {\bf R}^3$ be a helicoidal surface
given by \rf{tw} and assume that the surface has constant anisotropic mean curvature. If there exists an open set
$U\subset \Sigma$ such that $\nu_3\equiv 0$ on $U$, then all of
$X(\Sigma)$ is contained in a circular cylinder.
\end{lemma}
{\it Proof}  The Jacobi operator of a constant anisotropic mean curvature immersion is the elliptic self-adjoint operator
$$J[u]={\rm div}(D^2\gamma +\gamma I)\nabla u +\langle (D^2\gamma +\gamma I)d\nu, d\nu\rangle u\:.$$
In \cite{KP2005}, it is shown that $J[\nu_3]=0$ holds. By results of \cite{Hmndr},
the operator $J$ has the unique continuation property: if a solution vanishes identically on an open set, then the solution vanishes identically, so $\nu_3\equiv 0$ holds on all of $\Sigma$. We then see that $\eta_1\equiv 0$ holds and consequently 
$x^2+y^2\equiv $ constant so the surface is a cylinder. {\bf q.e.d.}\\[4mm]

{\it Proof of Theorem \rf{T1}} If  the surface is a helicoidal surface satisying \rf{H} 
there exists an open set on which $\nu_3\equiv 0$, then by the first lemma, the surface is contained in a round cylinder which is an example of a constant anisotropic mean curvature surface. 

Likewise if the surface is  helicoidal and there exists an open set on which $\nu_3\equiv 0$ holds, then by the second lemma, the surface is a vertical circular cylinder. For the cylinder $\eta_1\equiv 0$, so $\eta_2\equiv R \in {\bf R}$
and so \rf{H} holds with $-A=\Lambda R^2+R/\mu_2(0)$.
From now on, we  assume that no  open set exists on which $\nu_3$ vanishes.

The non parametric and twizzler representations of the surface and its normal give
\begin{equation}
\label{X}
X=(re^{i\theta}, g(r)+\lambda\theta)=( (x+iy)e^{-i\omega \vartheta}, \vartheta +c)\:\end{equation}
\begin{equation}
\label{nu}\nu=\frac{ (-(rg_r+i\lambda)e^{i\theta},r)}{\sqrt{r^2(1+g_r^2)+\lambda^2}}=\frac{(ie^{i(\phi-\omega\vartheta)}\:, \omega \eta_1)}{\sqrt{1+\omega^2\eta_1^2}}\:.\end{equation}
From these, we obtain the equalities
$$\frac{\nu_1+i\nu_2}{x_1+ix_2}=\frac{-(rg_r+i\lambda)}{r\sqrt{1+\omega^2\eta_1^2}}=\frac{ie^{i\phi}}
{(x+iy)\sqrt{1+\omega^2\eta_1^2}}\:.$$
Using the equation $x+iy=-(\eta_1+i\eta_2)e^{i\phi}$, we obtain
\begin{equation}
\label{q}
\frac{ -(rg_r+i\lambda)}{\sqrt{r^2(1+g_r^2)+\lambda^2}}= \frac{-(\eta_2+i\eta_1)}{(\eta_1^2+\eta_2^2)\sqrt{1+\omega^2\eta_1^2}}
\:.\end{equation}
By using the invariance of $X$ under the group action, we find $\lambda=-1/\omega$ and so from \rf{q}, we can conclude
$$rg_r=\frac{-\eta_2}{\omega \eta_1}\:.$$

It is clear from \rf{xiy} and \rf{X} that $r^2=x^2+y^2=\eta_1^2+\eta_2^2$ holds. Using
\rf{q} and the fact that $\nu_3=\omega\eta_1/\sqrt{1+\omega \eta_1^2}$, we see that \rf{tc}
is equivalent to \rf{int} with $2C=-A$.

 This verifies the conclusion of the theorem on any open set in the surface which can be represented as a graph over a horizontal plane, i.e. on any open set on which $\nu_3$ does not vanish.

The set $\{\nu_3=0\}$ is a closed set with empty interior. We write its compliment as 
$U_1\cup...\cup U_n$ where $U_i$ is open and connected. On each $U_i$, an equation of the form
\rf{H} holds where possible $A=A_i$ depends on $i$.  By considering a sequence of points, $p_k\in U_i$ with
$p_k\rightarrow p\in \partial U_i\cap \partial U_j$, we have $\eta_1(p_k)\rightarrow 0$ and
$-A_i= \Lambda \eta_2(p)+\frac{2\eta_2(p)}{\mu_2(0)}$. By considering a similar sequence in $A_j$, this shows that $A_i=A_j$ and so all of the 
$A_i's$ have a common value $A$.  {\bf q.e.d.}

\begin{figure}[h!]
\label{A}
\begin{minipage}[c]{.45\textwidth}
    \begin{center}
\includegraphics[width=40mm,height=55mm,angle=0]{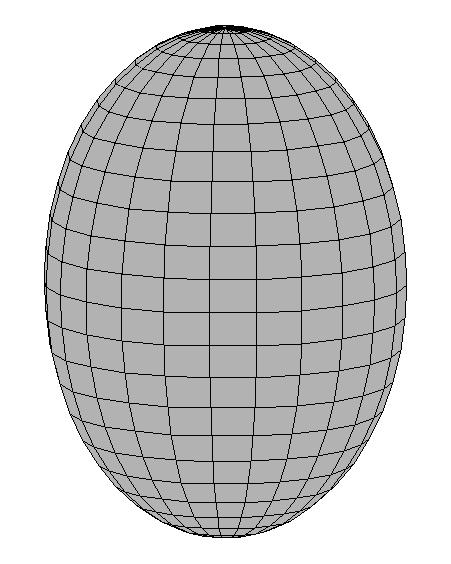}

% \label{c1}
    \end{center}
\end{minipage}
  %\hfill
  \begin{minipage}[c]{.45\textwidth}
    \begin{center}
\includegraphics[width=55mm,height=55mm,angle=0]{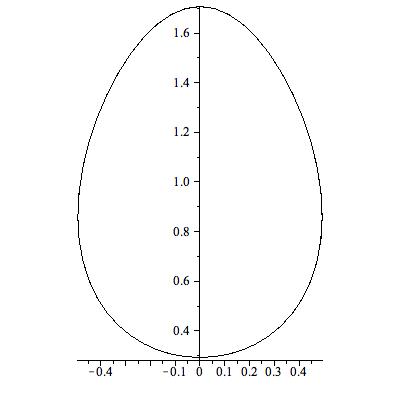}

\end{center}
 \end{minipage}

   \begin{minipage}[c]{.45\textwidth}
  \begin{center}
        \includegraphics[width=55mm,height=55mm,angle=0]{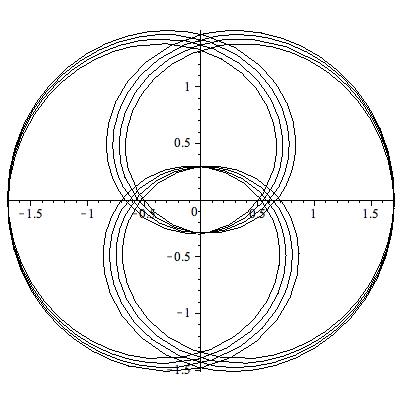}

     % \label{c2}
    \end{center}
   \end{minipage}
    \begin{minipage}[c]{.45\textwidth}
  \begin{center}
        \includegraphics[width=65mm,height=65mm,angle=0]{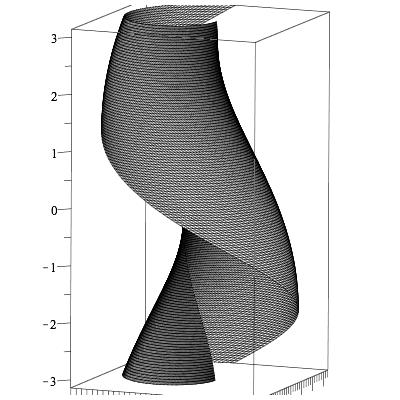}
 % \label{c2}
    \end{center}
   \end{minipage}
\caption{Upper left: Wulff shape for $\gamma=1+0.2\nu_3^2$. Upper right: treadmill sled
$(\Lambda,A,\omega)=(1,0.5,1)$. Bottom: Corresponding generating curve and twizzler surface.}
   \end{figure}
\begin{figure}[h!]
\label{B2}
\begin{minipage}[c]{.45\textwidth}
    \begin{center}
\includegraphics[width=40mm,height=55mm,angle=0]{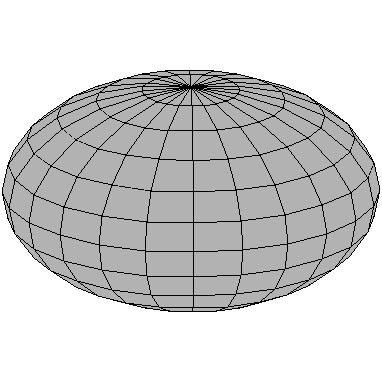}
%\caption{ Wulff shape for  $\gamma:=1-.03\nu_3^2$.}
% \label{c1}
    \end{center}
\end{minipage}
  %\hfill
  \begin{minipage}[c]{.45\textwidth}
    \begin{center}
\includegraphics[width=55mm,height=55mm,angle=0]{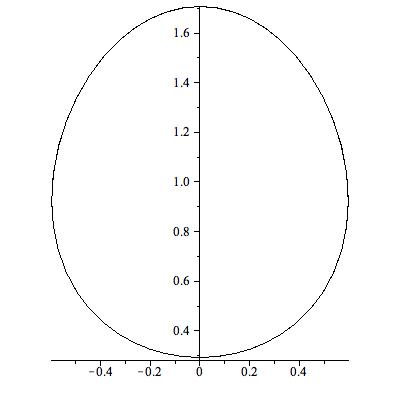}
%\caption{}
% \label{c2}
    \end{center}
   \end{minipage}
 % \caption{Wulff shape (left) and anisotropic catenoid (right)}
 % \label{dels1}
  %\end{figure} %\hfill
   \begin{minipage}[c]{.45\textwidth}
  \begin{center}
        \includegraphics[width=55mm,height=55mm,angle=0]{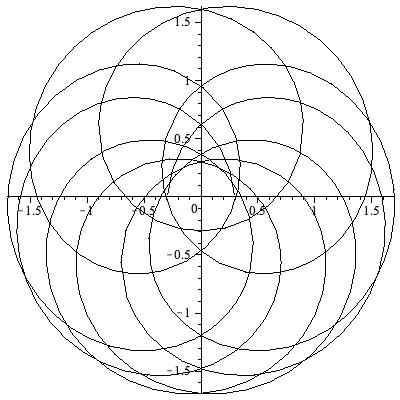}
%\caption{}

     % \label{c2}
    \end{center}
   \end{minipage}
    \begin{minipage}[c]{.45\textwidth}
  \begin{center}
        \includegraphics[width=70mm,height=70mm,angle=0]{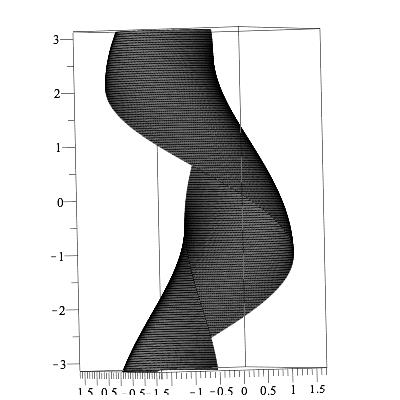}
 % \label{c2}
    \end{center}
   \end{minipage}
\caption{Upper left: Wulff shape for $\gamma=1-0.3\nu_3^2$. Upper right: treadmill sled
$(\Lambda,A,\omega)=(1,0.5,1)$. Bottom: Corresponding generating curve and twizzler surface.}
   \end{figure}
\section{Other anisotropies}

We will next discuss examples of  helicoid like zero mean curvature surfaces for other anisotropies. These will be based on the Wulff shapes:
$$W_p:=\{|\xi_1|^p+|\xi_2|^p+|\xi_3|^p=1\}\:,$$
where $p>1$. For $p>2$ (resp. $1<p<2$) $W$ is called a super (resp. sub) ellipsoid.  For convenience, we will take $p$ to be an even integer and as usual $q$ will denote the conjugate exponent defined by
$p^{-1}+q^{-1}=1$. 
The functional corresponding to the Wulff shape $W_p$ assigns to a surface having normal $\nu$ the value
$${\cal F}=\int_\Sigma |\nu|_q\:d\Sigma\:.$$
When $\Sigma$ is in non parametric form $z=z(x,y)$, we have
$${\cal F}=\int \bigl(1+ z_x^q+z_y^q\bigr)^{1/q}\:dxdy\:.$$
The Euler Lagrange equation for a  non-parametric equilibrium surface  $z=z(x,y)$  is given by
\begin{equation}
\label{ELq}
M_p[z]:= {\rm Div}_0\biggl([1+z_x^{q}+z_y^{q}]^\frac{1-q}{q}(z_x^{q-1}, z_y^{q-1} \biggr)=0\:
\end{equation}

Let $\omega=(|x|^p+|y|^p)^{1/p}$ .  We let $z=z(\omega)$ be a solution
of $M_p=0$. This will be the case when the graph of $z$ is a part of a (generalized) anisotropic catenoid for the Wulff shape $W_p$. We compute
$$z_x=z_\omega \omega_x=z_\omega \omega^{1-p}x^{p-1},\quad z_y=z_\omega \omega_x=z_\omega \omega^{1-p}y^{p-1}\:.$$
A straightforward computation shows that \rf{ELq} gives
\begin{equation}
\label{ELq2}
{\rm Div_0}\biggl( [\frac{z_\omega}{(1+z_\omega^q)^{1/q}}]^{q-1}\:\omega^{-1}\:(x,y)  \biggr)=0
\end{equation}
With the obvious notation, we will express this as
$${\rm Div}\biggl( f(\omega)(x,y)  \biggr)=0\:.$$

It follows that there is (locally away from $(0,0)$) a function $w(x,y)$ with $Dw=f(\omega)( -y,x)$, i.e.
\begin{equation}
\label{der}
w_x=-f(\omega)y\:,\qquad w_y=f(\omega)x\:.
\end{equation}
 We claim that
$$M_q[w]=0\:$$
holds, i.e. the graph of $w$ has zero anisotropic mean curvature for the functional whose Wulff shape
is $W_q$.  Of course, $W_q$ is just the unit sphere in the dual space to $({\bf R}^3,\:|\cdot|_p)$.

We have from \rf{der} that
$$[1+w_x^p+w_y^p]^\frac{1-p}{p}=[1+(\omega f(\omega))^p]^\frac{1-p}{p}$$
and (recall $q-1$ is odd!)
$$(w_x^{p-1},w_y^{p-1})=(f(\omega))^{p-1}(-y^{p-1}, x^{p-1})\:.$$
Combining these facts, we can write
$$M_q[w]={\rm Div_0}\biggl(\Gamma(\omega)(-y^{p-1}, x^{p-1}) \biggr)  \:.$$
for a suitable function $\Gamma$. Note that ${\rm Div}(-y^{q-1}, x^{q-1})=0$ and so
$$M_q[w]=\Gamma'(\omega)\langle \nabla\omega , (-y^{p-1}, x^{p-1})\rangle =\Gamma'(\omega)\omega^{p-1}
\langle (x^{p-1}, y^{p-1}) , (-y^{p-1}, x^{p-1})\rangle=0\:,$$
which proves the claim.

If we integrate using \rf{der}, we get
$$\oint_{\omega=a}dw=\oint_{\omega=a} f(\omega)(-ydx+xdy)=f(a)\oint_{\omega=a}
(-ydx+xdy)=2{\rm Area }(\{\omega \le a\})\ne 0\:,$$
so $dw$ has a non zero period and its graph, like that of a helicoid, is not single valued. 
Since $dw=  f(\omega)(-ydx+xdy)=0$ along radial lines, the graph of $w$ is a ruled surface. It is not difficult to compute $w$ explicitly
for a fixed value of $p$. We do not supply graphics of the resulting surfaces since they closely resemble the classical helicoid.

The functionals used above are a special case of a more general construction. Let 
$||\cdot ||$ be any smooth norm on ${\bf R}^3$. We take $W$ to be the unit sphere in this norm, i.e. $W=\{||\xi||=1\}$. The corresponding functional is 
$${\cal F}:=\int_\Sigma ||\nu||_*\:d\Sigma\:.$$
Here $||\cdot ||_*$ denotes the dual norm 
$$||x||_*=\sup_{||\xi||=1} \langle x,\xi\rangle\:.$$
Of course, this gives rise to a ``dual'' functional ${\cal F}^{\:*}$ whose Wulff shape is $W^*=\{||x||_*=1\}$.
We will always identify $({\bf R}^3)^*$ with $ {\bf R}^3$ by using the standard inner product.

We will now suppose that the norm $||\cdot ||$ has a special form. Let $|\cdot |_H$ and
$|\cdot |_V$ be smooth norms on ${\bf R}^2$ which we refer to as the horizontal and vertical norms
respectively. Then we assume that $||(a,b,c)||=|(|(a,b)|_H,c)|_V$. This norm has the property
that generic level sets of the height function of its unit sphere are all homothetic. It is not difficult to see that the dual norm with have the same form and that $||(A,B,C)||_*=|( |(A,B)|_{H_*},C)|_{V_*}$.

We will derive the Euler-Lagrange equation for the functional ${\cal F}$ for a surface in non-parametric form $w=w(x,y)$. Because of the homogeneity of the norm, we get
\begin{eqnarray*}
{\cal F}^{\*}=\int_\Sigma ||\nu||\: d\Sigma &=& \int_G ||(-\nabla w,1)||\:d^2x \\
                                                              &=& \int_G |(|\nabla w|_H,1)|_V\:d^2x\:.
\end{eqnarray*}
For $ p\in {\bf R}^2$, set $\Psi(p):=|p|_H$, $\Phi(p):=|p|_V$. Then
$$\delta \bigl[\Phi((\Psi(\nabla w),1))\bigr]=\nabla \Phi_{(\Psi(\nabla w),1)}\cdot \bigl( (\nabla \Psi |_{\nabla w}\cdot \nabla {\dot w}),0\bigr)= \bigr[ \nabla \Phi_{(\Psi(\nabla w),1)}\cdot (1,0)\bigr] (\nabla \Psi |_{\nabla w}\cdot \nabla {\dot w})           \:.$$
It follows that for a variation of $w$, $w\rightarrow w+\epsilon {\dot w}$, 
\begin{eqnarray*}
\delta {\cal F}^* &=&\int_G \bigr[ \nabla \Phi_{(\Psi(\nabla w),1)}\cdot (1,0)\bigr] (\nabla \Psi |_{\nabla w}\cdot \nabla {\dot w}) \:d^2x\\
&=& -\int_G {\dot w}\: {\rm Div}_0 \biggl(\bigr[ \nabla \Phi_{(\Psi(\nabla w),1)}\cdot (1,0)\bigr] \nabla \Psi |_{\nabla w} \biggr)\:d^2x.
\end{eqnarray*}
So the Euler-Lagrange equation
\begin{equation}
\label{ELN}  {\rm Div}_0 \biggl(\bigr[ \nabla \Phi_{(\Psi(\nabla w),1)}\cdot (1,0)\bigr] \nabla \Psi |_{\nabla w} \biggr)=0\:,\end{equation}
expresses the vanishing of the anisotropic mean curvature of the graph for the functional having Wulff shape given by $W^*:=\{||X||_*=1\}$. 

Let $\Psi^*(p):=|p|_{H*}$,  $\Phi^*(p):=|p|_{V*}$.  We seek a solution of \rf{ELN} of the form $w=w(\Psi^*(x,y))$. We have
\begin{equation}
\label{giga1}
\nabla \Psi |_{\nabla w}=\nabla \Psi |_{w'(\Psi^* )\nabla \Psi^*(x,y)}=\nabla \Psi |_{\nabla \Psi^*(x,y)}= \frac{(x,y)}{\Psi^*(x,y)}\:.\end{equation}
The second equality above follows from the homogeneity of $\Psi^*$ while the third equality follows from equation (1.7.9) of \cite{G}. Also, we have $\Psi(\nabla w)=|w'(\Psi^*)|\Psi(\nabla
\Psi^*)=|w'(\Psi^*)|$ by (1.7.8) of \cite{G}. We then obtain  that, with
$g(\Psi^*):=[\nabla \Phi^*|_{(w'(\Psi^*),1)}\cdot (1,0)] /\Psi^*$,  $w$ must satisfy

\begin{equation}
\label{ELN2}  {\rm Div}_0\bigl(g(\Psi^*)(x,y)\bigr)=0\:.
\end{equation}

Any solution of this equation has zero anisotropic mean curvature for the functional with Wulff shape $W^*$ and the solution has cross-sections which are rescalings of generic cross sections
of $W^*$. It is therefore called an {\it anisotropic catenoid}. These surfaces were first constructed in \cite{KP2008} by another method.

We will use the equation \rf{ELN2} to construct a ``helicoid'' for the functional ${\cal F}$. 
From \rf{ELN2}, we have, away from the origin, the local existence of a function $\alpha$ satisfying 
\begin{equation}
\label{dual}
\nabla \alpha=g(\Psi^*)J(x,y):=g(\Psi^*)(-y,x)\:.
\end{equation}
\begin{theorem}
Assume that  
\begin{equation}
\label{J}
\Psi^*(Jv)=\Psi^*(v)\:,
\end{equation}
holds for all $v\in {\bf R}^2$. Then, away form the origin, $\alpha$ satisfies the dual equation
\begin{equation}
\label{ELN3}  {\rm Div}_0 \biggl(\bigr[ \nabla \Phi^*_{(\Psi^*(\nabla \alpha),1)}\cdot (1,0)\bigr] \nabla \Psi^* |_{\nabla \alpha} \biggr)=0\:,\end{equation}
so the graph of $\alpha$ has zero anisotropic mean curvature for the functional with Wulff shape
$W$.  Furthermore, $\alpha$ is multivalued in any punctured neighborhood of the origin and the graph of $\alpha$  is ruled by horizontal lines.
\end{theorem}
{\it Proof.} First we have 
$$\nabla \Psi^*_{\nabla \alpha}=\nabla \Psi^*_{g(\Psi^*)J(x,y)}=
J\nabla \Psi^*_{(x,y)}\,$$
 by \rf{J} and the fact that $\nabla \Psi^*$ is homogeneous of degree zero.
Also, we have
$$\Psi^*(\nabla \alpha)=|g(\Psi^*)|\Psi^*(J(x,y))=|g(\Psi^*)|\Psi^*((x,y))=:h(\Psi^*)\:,$$
again using \rf{J}. We can therefore define a function:
$$\eta(\Psi^*):=[ \nabla \Phi^*_{(\Psi^*(\nabla \alpha),1)}\cdot (1,0)]
=\nabla \Phi^*_{(h(\Psi^*),1)}\cdot (1,0)]\:,$$
from which we get
$${\rm Div}_0 \biggl(\bigr[ \nabla \Phi^*_{(\Psi^*(\nabla \alpha),1)}\cdot (1,0)\bigr] \nabla \Psi^* |_{\nabla \alpha} \biggr)={\rm Div}_0\biggl(\eta(\Psi^*)J\nabla \Psi^*\biggr)=0\:,$$
since ${\rm Div}_0 J\nabla f=0$ for any smooth function $f$.
Note that for any positive constant $c$,
$$\int_{\Psi^*=c} d\alpha =g(c) 2{\rm Area}(\Psi^*\le c)\:.$$
It is easy to see that $g$ is not identically zero, so $\alpha$ is multivalued.
The final statement of the theorem follows from
$$\nabla \alpha \cdot (x,y)=g(\Psi^*)(-y,x)\cdot (x,y)\equiv 0,$$
so the height function is constant on radial lines. {\bf q.e.d}

\begin{flushleft}
Chad  K{\footnotesize UHNS}\\
Department of Mathematics\\
Montgomery College\\
Rockville, MD, 20850\\
U.S.A.\\
E-mail: 
Chad.Kuhns@montgomerycollege.edu
\end{flushleft}

\begin{flushleft}
Bennett P{\footnotesize ALMER}\\
Department of Mathematics\\
Idaho State University\\
Pocatello, ID 83209\\
U.S.A.\\
E-mail: palmbenn@isu.edu
\end{flushleft}

\end{document}